\documentclass[a4paper,12pt,leqno]{amsart}

\usepackage[latin1]{inputenc}
\usepackage{xspace}

\newcommand{\R}{\mathbb{R}}
\newcommand{\PP}{\mathbb{P}}
\newcommand{\E}{\mathbb{E}}

\newcommand{\var}{\mathop{\rm Var}}

\newcommand{\unsur}[1]{{\frac{1}{#1}}}
\def\un#1{{{\,\mathbf{1}}_{({#1})}}}
\def\valabs#1{{\left\vert {#1} \right \vert}}
\def\esp#1{{{\E}\left [ {#1} \right ]}}
\def\etp#1{{\left ( {#1} \right )}}
\def\etc#1{{\left [ {#1} \right ]}}
\def\crochet#1{{\left < {#1} \right >}}
\def\esperance#1#2{{{\E}_{#1}\left [ {#2} \right ]}}
\def\prob#1{{{\PP}\left ( {#1} \right )}}

\def\ens#1{{\left\{#1\right\}}}

\def\dessus#1#2{\mathord{\mathop{\kern 0pt #2}\limits^#1}}

\newcommand{\bit}{\begin{itemize}}
\newcommand{\eit}{\end{itemize}}
\newcommand{\ben}{\begin{enumerate}}
\newcommand{\een}{\end{enumerate}}

\newcounter{moncompteur}

\newenvironment{myenumerate}%
{\begin{list}{\arabic{moncompteur}. }{\usecounter{moncompteur}%
\setlength{\leftmargin}{0pt}%
\setlength{\labelwidth}{0pt}%
\setlength{\listparindent}{0pt}%
\setlength{\labelsep}{0pt}}}%
{\end{list}}

\def\bmen{\begin{myenumerate}}
\def\emen{\end{myenumerate}}

\newcommand{\intof}{{\int_0^\infty}}

\newcommand{\intot}{{\int_0^t}}

\newcommand{\undemi}{\frac{1}{2}}

\newcommand{\Frond}{{\mathcal F}}
\newcommand{\Grond}{{\mathcal G}}

\newcommand{\Nrond}{{\mathcal N}}

\newtheorem{theorem}{Theorem}

\newtheorem{proposition}[theorem]{Proposition}

\theoremstyle{definition}

\theoremstyle{remark}

\newcommand{\as}{almost surely\xspace}

\newcommand{\gesp}[1]{\mathbf{E}\etc{#1}}
\newcommand{\gprob}[1]{\mathbf{P}\etp{#1}}

\begin{document}
\title[Strong disorder implies strong localization for DPRE]{  Strong disorder implies strong localization for directed polymers in a random environment}

\author{Philippe Carmona}
\address{Philippe Carmona\\
Laboratoire Jean Leray, UMR 6629,
Universit{\'e} de Nantes, BP 92208\\
F-44322 Nantes Cedex 03
BP 
}
\email{philippe.carmona@math.univ-nantes.fr}
\author{Yueyun Hu}
\address{Yueyun Hu\\
D\'epartement de Math\'ematiques\\
 (Institut Galil\'ee, L.A.G.A. UMR 7539)\\
Universit\'e Paris 13
99 Av. J-B Cl\'ement, 93430 Villetaneuse, France}
\email{yueyun@math.univ-paris13.fr}

\begin{abstract}
In this note we show that in any dimension $d$, the strong disorder property implies the strong localization property. This is established for a continuous time model of directed polymers in a random environment : the parabolic Anderson Model.
\end{abstract}


%\thanks{}

\date{\today}\vfuzz=4pt \hfuzz=40pt
\maketitle

\section{Introduction}
Let $\omega=(\omega(t))_{t\ge 0}$ be the simple continuous time random walk  on the $d$-dimensional lattice $\mathbb{Z}^d$, with jump rate $\kappa >0$, defined on a probability space $(\Omega,\Frond,\PP)$. We consider an \emph{ environment } $B=(B_x(t),t\ge 0, x \in\mathbb{Z}^d)$ made of independent standard Brownian motions $B_x$ defined on another probability space $(H,\Grond,\mathbf{P})$.

For any $t>0$ the (random) \emph{polymer measure} $\mu_t$ is the probability defined on the path space $(\Omega,\Frond)$ by
\begin{equation*}
  \label{eq:1}
  \mu_t(d\omega)= \unsur{Z_t} e^{\beta H_t(\omega) - t \beta^2/2} \PP(d\omega)
\end{equation*}
where $\beta \ge 0$ is the inverse temperature, the Hamiltonian is
\begin{equation*}
  \label{eq:2}
  H_t(\omega) = \intot dB_{\omega(s)}(s)
\end{equation*}
and the \emph{partition function} is 
\begin{equation*}
  \label{eq:3}
  Z_t = Z_t(\beta) = \esp{e^{\beta H_t(\omega) - t \beta^2/2}}\,,
\end{equation*}

where $\esp{}$ denotes expectation with respect to $\PP$.

\smallskip

Erwin Bolthausen~\cite{MR1006293} was the first to establish that $(Z_t)_{t\ge 0}$ was a positive martingale, converging \as to a finite random variable $Z_\infty$, satisfying a {\emph zero-one} law : $\gprob{Z_\infty >0}\in\ens{0,1}$. We shall say that there is \emph{strong disorder} if $Z_\infty=0$ \as, and \emph{weak disorder} if $Z_\infty >0$ \as.

Another martingale argument, based on a supermartingale decomposition of $\log Z_t$, enabled Carmona-Hu~\cite{MR2073415}, then Comets-Shiga-Yoshida~\cite{MR1996276,MR2073332}, and Rovira-Tindel \cite{MR2129770}, to show the equivalence between strong disorder and \emph{weak-localization} : 
\begin{equation*}
  \label{eq:4}
  Z_\infty =0\; a.s. \quad \iff \intof \mu_t^{\otimes 2}(\omega_1(t)=\omega_2(t))\, dt = +\infty \quad a.s. \,,
\end{equation*}
where $\omega_1,\omega_2$ are two independent copies of the random walk $\omega$, considered under the product polymer measure $\mu_t^{\otimes 2}$:
\begin{equation*}
  \label{eq:5}
  \mu_t^{\otimes 2}(d\omega_1,d\omega_2) = \unsur{Z_t^2} e^{\beta(H_t(\omega_1) + H_t(\omega_2)) - t \beta^2}\; \PP^{\otimes 2}(d\omega_1,d\omega_2)\,.
\end{equation*}

Let us define \emph{strong localization} as the existence of a constant $c>0$ such that 
\begin{equation*}
  \label{eq:6}
  \limsup_{t\to +\infty}\,  \sup_x \mu_t(\omega(t)=x) \ge c \quad a.s.
\end{equation*}
This property implies the existence of highly favored sites, in contrast to the simple random walk ($\beta=0$) for which $\sup_x \prob{X_t=x} \sim C t^{-d/2} \to 0$.
Carmona-Hu~\cite{MR2073415}, and then Comets-Shiga-Yoshida~\cite{MR2073332}, showed that in dimension $d=1,2$, for any $\beta>0$, there was not only strong disorder but also strong localization.

\bigskip
We shall prove in this note the

\begin{theorem}
In any dimension $d$,  strong disorder implies strong localization.
\end{theorem}

For sake of completeness, let us state yet another localization property. The \emph{free energy} is the limit
\begin{equation*}
  \label{eq:7}
  p(\beta) = \lim_{t\to +\infty} \unsur{t} \log Z_t\,,
\end{equation*}
where the limit can be shown to hold \as and in  every $L^p$, $p\ge 1$ (see e.g. ~\cite{MR2073332}). The  function $p(\beta)$ is continuous, non increasing on $[0,+\infty[$,  $p(\beta)\le 0$, $p(0)=0$, so there exists a critical inverse temperature $\beta_c \in [0,+\infty]$ such that:
\begin{equation*}
  \label{eq:8}
  \begin{cases}
    p(\beta)=0 & \text{if } 0 \le \beta \le \beta_c\,; \\
p(\beta)< 0 & \text{if } \beta > \beta_c\,.
  \end{cases}
\end{equation*}
When $p(\beta)<0$ we say that the system has the \emph{very strong disorder} property. We shall prove that (see equation~\eqref{par:eq:7}):
\begin{equation*}
  p(\beta) = - \frac{\beta^2}{2} \lim_{t \to +\infty} \unsur{t} \intot \mu_s^{\otimes 2} (\omega_1(s)=\omega_2(s))\,ds \quad a.s.
\end{equation*}
Therefore  there is \emph{very strong disorder} if and only if there exists a constant $c>0$ such that almost surely:
\begin{equation*}
   \lim_{t \to +\infty} \unsur{t} \intot \mu_s^{\otimes 2} (\omega_1(s)=\omega_2(s))\,ds =c
\end{equation*}

 The recent beautiful result of Comets-Vargas~\cite{math.PR/0510525}, that is $\beta_c=0$ in dimension $d=1$, strengthen our belief in the 
 \begin{equation*}
   \label{eq:10}
   \text{\texttt{Conjecture} : } \quad \text{very strong disorder} \iff \text{strong disorder}
 \end{equation*}
Proving this conjecture would unify all these notions of disorder and localization.

\medskip

Eventually, let us end this rather lengthy introduction by making clearer the connection with the parabolic Anderson model (see Carmona and Molchanov~\cite{MR1185878} or Cranston, Mountford and Shiga~\cite{MR1980378}). The point to point partition functions
\begin{equation*}
  \label{eq:11}
  Z_t(x,y) = \esperance{x}{e^{\beta H_t(\omega) - t \beta^2/2} \un{\omega(t)=y}}
\end{equation*}
satisfy the stochastic partial differential equation (see Section~2) 
\begin{equation*}
  \label{eq:12}
  dZ_t(0,x) = L Z_t(0,.)(x) \, dt + \beta Z_t(0,x) \, dB_x(t)\,,
\end{equation*}
where $L=\kappa \Delta$ is the generator of the simple random walk $\omega$ with jump rate $\kappa$, that is $\Delta$ is the discrete Laplacian.

\bigskip
Let us explain now the structure of this paper. Section 2 is devoted to the study of the partition function as a martingale, and we prove that its asymptotics are governed by the asymptotics of the overlap $I_t=\mu_t^{\otimes 2}(\omega_1(t)=\omega_2(t))$.

An important fact is that $I_t$ itself is a semimartingale. In Section 3 we establish a decomposition of $I_t$ which is not its canonical semimartingale decomposition (this decomposition can be obtained via the parabolic Anderson equation\eqref{eq:11}). In fact this decomposition looks a lot like a renewal equation involving the overlap for the simple random walk : it is the basic ingredient of our proof of the main result, since it is in this decomposition that we inject our knowledge of the behaviour of the overlap for simple random walk.
\section{The partition function}

Without loss in generality we can work on the canonical path space $\Omega$ made of $\omega:\R^+ \to \mathbb{Z}^d$, càdlàg, with a finite number of jumps in each finite interval $\etc{0,t}$. We endow $\Omega$ with the canonical sigma-field $\Frond$ and the family of laws $(\PP_x,x\in\mathbb{Z}^d)$ such that under $\PP_x$, $(\omega(t))_{t\ge 0}$ is the simple random walk starting from $x$, with generator $L=\kappa \Delta$. With these notations, we consider, attached to each path $\omega \in \Omega$, the exponential martingale
\begin{equation*}
  \label{par:eq:1}
  M^\omega_t = \exp(\beta H_t(\omega) - t \beta^2/2) = 1 + \beta \intot M^\omega_s \, dB_{\omega(s)}(s)\,,
\end{equation*}
with respect to the filtration $\Grond_t = \sigma(B_x(s), s\le t, x\in\mathbb{Z}^d)$.
We have $Z_t=\esp{M^\omega_t}$ and thus the
\begin{proposition}
The process $(Z_t)_{t\ge 0}$ is a continuous positive $\Grond_t$ martingale with quadratic variation
\begin{equation*}
  \label{par:eq:2}
  d\crochet{Z,Z}_t = Z_t^2 \beta^2 \, I_t\, dt\,,\quad\text{with}\quad I_t = \mu_t^{\otimes 2}(\omega_1(t)=\omega_2(t))\,.
\end{equation*}
  
\end{proposition}
\begin{proof}
  We know that linear combinations of martingales are martingales. This extends easily to probability mixtures of martingales. Indeed, let $0\le s\le t$ and let $U$ be positive bounded and $\Grond_s$-measurable. Then, by Fubini-Tonelli's theorem :
  \begin{align*}
    \gesp{Z_t U} &= \gesp{\esp{ M^\omega_t} U} = \esp{\gesp{M^\omega_t U}} \\
&= \esp{\gesp{M^\omega_s U}} &\text{($M^\omega$ is a martingale)}\\
&= \gesp{\esp{M^\omega_s} U} = \gesp{Z_s U}\,.
  \end{align*}
Observe that if $\omega_1,\omega_2$ are paths, then we can compute the quadratic covariation
\begin{equation*}
  \label{par:eq:3}
d\crochet{M^{\omega_1},M^{\omega_2}}_t = M^{\omega_1}_t M^{\omega_2}_t \beta^2 \un{\omega_1(t)=\omega_2(t)}\, dt.
\end{equation*}
Therefore, we have formally:
\begin{align*}
  d\crochet{Z,Z}_t &= d\crochet{\int \PP(d\omega_1) M^{\omega_1}, \int \PP(d\omega_2) M^{\omega_2}}_t \\
&= \int \PP^{\otimes 2}(d\omega_1,d\omega_2)d\crochet{M^{\omega_1},M^{\omega_2}}_t  \\
&= \beta^2 Z_t^2  \unsur{Z_t^2} \int \PP^{\otimes 2}(d\omega_1,d\omega_2)M^{\omega_1}_t M^{\omega_2}_t  \un{\omega_1(t)=\omega_2(t)}\, dt \\
&=  Z_t^2 \beta^2 \, I_t\, dt.
\end{align*}
This again can be made rigorous by writing $N_t = Z_t^2 -\beta^2 \intot Z_s^2 I_s \, ds$ as a probability mixture of martingales:
\begin{equation*}
  \label{par:eq:4}
  N_t =  \int \PP^{\otimes 2}(d\omega_1,d\omega_2) (M^{\omega_1}_t M^{\omega_2}_t - \beta^2 \intot M^{\omega_1}_s M^{\omega_2}_s \un{\omega_1(s)=\omega_2(s)}\, ds)\,.
\end{equation*}
\end{proof}
The positive martingale $Z_t$ converges \as to a positive finite random variable $Z_\infty$. We refer to any of \cite{MR1006293,MR2073332,MR1939654} for a proof of the following zero-one law.
\begin{proposition}
$$  \gprob{Z_\infty =0} \in \ens{0,1}\,.$$
\end{proposition}

We can now show the equivalence between strong disorder and weak localization.
\begin{proposition}
  The supermartingale $\log Z_t $ has the decomposition
  \begin{equation*}
    \label{par:eq:5}
\log Z_t = M_t - \undemi A_t
  \end{equation*}
with $(M_t)_{t\ge 0}$ a continuous martingale of quadratic variation
\begin{equation*}
  \label{par:eq:6}
\crochet{M,M}_t = A_t = \beta^2 \intot I_s\, ds\,.
\end{equation*}
Consequently:
\begin{itemize}
\item either $Z_\infty =0$ and $\intof I_s \, ds = +\infty$ \as ;
\item or $Z_\infty >0$ and $ \intof I_s \, ds < +\infty$ \as.
\end{itemize}

In both cases the free energy is given by
\begin{equation}
  \label{par:eq:7}
p(\beta) = - \frac{\beta^2}{2} \lim_{t\to +\infty} \unsur{t} \intot I_s\, ds =  - \frac{\beta^2}{2} \lim_{t\to +\infty} \unsur{t} \intot \gesp{I_s}\, ds\,.
\end{equation}
\end{proposition}

\begin{proof} One can even prove (see~\cite{MR1939654}) that weak disorder is equivalent to the uniform integrability of the martingale $(Z_t)_{t\ge 0}$.
\medskip

Itô's formula yields :
\begin{equation*}
  \label{par:eq:8}
\log Z_t = \intot \frac{d Z_s}{Z_s} - \undemi \intot \frac{d\crochet{Z,Z}_s}{Z_s^2} = M_t - \undemi \beta^2 \intot I_s\, ds = M_t - \undemi A_t. 
\end{equation*} 

Therefore,
\begin{itemize}
\item On $\ens{A_\infty=\crochet{M,M}_\infty<+\infty}$ the martingale $M_t$ converges \as, $M_t \to M_\infty$ so $\log Z_t \to M_\infty - \undemi A_\infty$ and $Z_\infty >0$ \as, and $p(\beta)=\lim_{t\to +\infty} \unsur{t} \log Z_t = 0$.
\item On $\ens{A_\infty=\crochet{M,M}_\infty=+\infty}$, we have \as \, $\frac{M_t}{\crochet{M,M}_t} \to 0$ so $\frac{\log Z_t}{A_t} \to -\undemi$ and $ \log Z_t \to -\infty$, so $Z_\infty =0$. Furthermore, $p(\beta)=\lim_{t\to +\infty} \unsur{t} \log Z_t = -\undemi \lim_{t\to +\infty} \unsur{t} A_t$.
\end{itemize}

We conclude this proof by taking expectations:
$$ p(\beta) = \lim_{t\to +\infty} \unsur{t} \gesp{\log Z_t} = -\undemi \lim_{t\to +\infty} \unsur{t} \gesp{A_t} =   - \frac{\beta^2}{2} \lim_{t\to +\infty} \unsur{t} \intot \gesp{I_s}\, ds\,.$$

\end{proof}

The connection with the parabolic Anderson model is contained in the 
\begin{proposition}
  The point to point partition functions $(Z_t(0,x), t\ge 0, x \in \mathbb{Z}^d)$ 

satisfy the stochastic partial differential equation
\begin{equation*}
  dZ_t(0,x) = L Z_t(0,.)(x) \, dt + \beta\,Z_t(0,x) \, dB_x(t)\,,
\end{equation*}
where $L=\kappa \Delta$ is the generator of the simple random walk with jump rate $\kappa$, that is $\Delta$ is the discrete Laplacian.
\end{proposition}
\begin{proof}
  Let $p_t(x)=\prob{X_t=x}$ be the probability function at time $t$ of simple random walk. By Fubini's stochastic theorem and Markov property:
  \begin{align*}
    Z_t(0,x) &= \int \PP(d\omega) M^\omega_t \un{\omega(t)=x} \\
&= \int \PP(d\omega) \un{\omega(t)=x}(1+ \beta \intot M^\omega_s dB_{\omega(s)}(s)) \\
&= p_t(x) + \beta \intot \int \PP(d\omega) \un{\omega(t)=x} M^\omega_s dB_{\omega(s)}(s) \\
&=  p_t(x) + \beta \intot \int \PP(d\omega) p_{t-s}(\omega(s)-x) M^\omega_s dB_{\omega(s)}(s) \\
&= p_t(x) + \beta \intot Z_s \mu_s(p_{t-s}(\omega(s)-x)  dB_{\omega(s)}(s))\,.
  \end{align*}
We conclude by differentiating with respect to $t$, taking into account that
$$ \frac{d}{dt} p_t(x) = Lp_t(x)$$
In other words, we combine
$$ p_{t-s}(y) = \un{y=0} + \int_s^t Lp_{u-s}(y)\, du$$
and Fubini's stochastic theorem. (This result is just Feynman-Kac formula combined with time reversal of the continuous time random walk).
\end{proof}

\section{Itô's formula for the polymer measure}
\newcommand{\pn}{P^{\otimes n}}
\newcommand{\mom}{\mbox{\boldmath$\omega$}}

Let $(\pn_t)_{t\ge 0}$ be the semi-group of the Markov process $\mom(t)=(\omega_1(t), \ldots,\omega_n(t))$ constructed from $n$ independent copies of the simple random walk $(\omega(t))_{t\ge 0}$: 
if $f:\R^n \to \R$ is a bounded Borel function, then
$$ \pn_t f( x_1, \ldots, x_n) = \mathbb{E}_{x_1,\ldots,x_n}\etc{f(\omega_1(t),\ldots,\omega_n(t))}\,.$$

\begin{theorem}\label{thmitopoly}
  Let $f:\R^n \to \R$ be a bounded Borel function, and $t\ge t_0\ge 0$. Then,
  \begin{multline*}
    \mu^{\otimes n}_t\etc{f(\mom(t))}= \mu^{\otimes n}_{t_0}\etc{\pn_{t-t_0}f(\mom(t_0))}  \\
 + \beta^2 \sum_{i<j}\int_{t_0}^t \mu_s^{\otimes n}\etc{  \un{\omega_i(s)=\omega_j(s)} \pn_{t-s} f(\mom(s))}\, ds \\
 -n \beta^2 \sum_i \int_{t_0}^t \mu_s^{\otimes (n+1)}\etc{  \un{\gamma(s)=\omega_i(s)}  \pn_{t-s} f(\mom(s))}\, ds \\
 + \frac{n(n+1)}{2} \beta^2 \int_{t_0}^t \mu_s^{\otimes n}\etc{ \pn_{t-s} f(\mom(s))}\, I_s \, ds \\
+ \int_{t_0}^t \mu_s^{\otimes n}\etc{\pn_{t-s} f(\mom(s)) (\beta \sum_i  dB_{\omega_i(s)}(s) - n \frac{dZ_s}{Z_s})}\,,
  \end{multline*}
where $\gamma $ is an  extra independent copy of $\omega$. 
\end{theorem}

\begin{proof}
  Given paths $\omega_1,\ldots, \omega_n$, we let 
$$U_t = U_t(\omega_1,\ldots,\omega_n) = \frac{M^{\omega_1}_t \ldots M^{\omega_n}_t}{Z_t^n}\,.$$
We use the following easy computations of quadratic variations:
\begin{gather*}d\crochet{M^\gamma,M^\tau}_t = M^\gamma_t M^\tau_t \beta^2 \un{\gamma(t)=\tau(t)}\, dt \\
 d\crochet{M^\gamma,Z}_t = \beta^2 M^\gamma_t Z_t \mu_t\etc{\un{\omega(t) =\gamma(t)}}\,dt\,,\qquad
d\crochet{Z,Z}_t = Z_t^2\,\beta^2  I_t dt\,,\quad 
\end{gather*}

The classical Itô's formula yields:

\begin{align*}
U_t &= U_{t_0} + \int_{t_0}^t U_s \etp{\sum_{i=1}^n \beta  dB_{\omega_i(s)}(s) - n  \frac{dZ_s}{Z_s}} \\ 
 &+ \beta^2 \int_{t_0}^t U_s \etp{\sum_{i<j} \un{\omega_i(s)=\omega_j(s)} -n \sum_i \mu_s\etc{\un{\gamma(s)=\omega_i(s)}} +  \frac{n(n+1)}{2}I_s} ds\,,
\end{align*}
where in the last line $\mu_s$ acts on the generic path $\gamma$.
Since,
$$ \mu^{\otimes n}_t\etc{f(\mom(t))} = \int  f(\mom(t)) U_t(\mom)\, d\PP^{\otimes n}(\mom)$$
we conclude this proof by applying Fubini's theorem and Markov's property. For example,

\begin{align*}
  \int f(\mom(t))  U_{t_0}(\mom)\,d\PP^{\otimes n}(\mom) &=
\esp{f(\mom(t)) \frac{M^{\omega_1}_{t_0} \ldots M^{\omega_n}_{t_0}}{Z(t_0)^n}} \\
&= \unsur{Z(t_0)^n} \esp{\pn_{t-t_0}f(\mom(t_0)) M^{\omega_1}_{t_0} \ldots M^{\omega_n}_{t_0}} \\
&= \mu_{t_0}^{\otimes n}\etc{\pn_{t-t_0}f(\mom(t_0))}\,.
\end{align*}
\end{proof}

\section{Proof of the main result}

\newcommand{\espd}[1]{\mathbb{E}^{\otimes 2}\etc{#1}}
\newcommand{\probd}[1]{\mathbb{P}^{\otimes 2}\etp{#1}}

We assume that there is strong disorder so \as, $Z_\infty =0$ and $\intof I_s\, ds =+\infty$, and we shall show that for a certain $c_0>0$, $\limsup_{t \to +\infty} V_t \ge c_0 $ \as, with $V_t = \sup_x \mu_t(\omega(t)=x)$.

\medskip

Let $r(t)= \probd{\omega_1(t)=\omega_2(t)}$  and $R(t)=\intot r(s)\, ds$. In dimension $d=1,2$, $R(\infty)=+\infty$ so certainly $\beta^2 R(\infty) >1$. In dimension $d\ge 3$, $R(\infty)<+\infty$ and  Markov's property implies that $L_\infty = \intof \un{\omega_1(s)=\omega_2(s)}\, ds $ is under $\PP^{\otimes 2}$ an exponential random variable of expectation $R(\infty)$. Since, by Fubini's theorem,
\begin{align*}
  \gesp{Z_t^2} &= \espd{\gesp{e^{\beta(H_t(\omega_1) + H_t(\omega_2)) - t \beta^2}}}
\\&= \espd{ e^{\frac{\beta^2}{2} \var(H_t(\omega_1) + H_t(\omega_2)) -t \beta^2}} \\&=\espd{ e^{\beta^2 \intot \un{\omega_1(s)=\omega_2(s)}\, ds}}\,,
\end{align*}
the second moment method yields that if $\beta^2 R(\infty)<1$, then $sup_t \gesp{Z_t^2} =\espd{e^{\beta^2 L_\infty}} < +\infty$, so $Z_t$ is an $L^2$ bounded martingale, hence $\esp{Z_\infty}=1$ and $Z_\infty >0$ \as. Birkner~\cite{MR2041302} improved this result by using a conditional moment method : if $R(\infty)<+\infty$, then there exists $\beta^{-}_c > \unsur{\sqrt{R(\infty)}}$ such that for $\beta< \beta_c^{-}, Z_\infty >0$ \as. Hence, since we assumed strong disorder, we certainly have $\beta^2 R(\infty) >1$.

\medskip

Observe that since $V_t = \sup_x U_t(x)$ with $U_t(x)=\mu_t(\omega(t)=x)$, we have
\begin{align*}
 I_t &= \mu_t^{\otimes 2}(\omega_1(t)=\omega_2(t)) = \sum_x \mu_t^{\otimes 2}(\omega_1(t)=x=\omega_2(t))\\
&= \sum_x U_t(x)^2 \le V_t \sum_x U_t(x) = V_t
\end{align*}

and $I_t \ge V_t^2$. Therefore we shall show that \as, $\limsup_{t \to +\infty} I_t \ge c_0$. It is sufficient to prove that if $J_t= I_t \un{I_t \ge c_0}$ then for a constant $c_1>0$, 
\begin{equation*}
  \label{ito:eq:1}
  \limsup_{t\to +\infty} \frac{\intot J_s\, ds}{\intot I_s\, ds}\ge c_1 \quad \text{\as},
\end{equation*}
(indeed recall that $\intof I_s\, ds =+\infty$ \as).

We now have to choose $c_0>0$. Since $\beta^2 R(\infty)>1$, there exists $\epsilon_0 \in (0, \unsur{16})$ and $t_0>0$ such that $\beta^2 R(t_0) (1 - 4 \sqrt{\epsilon_0}) >1$. We let $c_0 = \epsilon_0 \inf_{0\le t \le t_0} r(t)$.

Let us apply now Itô's formula of Theorem~\ref{thmitopoly}, between $t-t_0$ and $t$, to the function $f(x_1,x_2)= \un{x_1=x_2}$:
\begin{align}\label{eq:14}
  I_t &=\mu_t^{\otimes 2}(f(\mom(t))) = N_{t_0,t} + \mu^{\otimes 2}_{t-t_0} \etc{P^{\otimes 2}_{t_0} f(\mom(t-t_0))}  \\
&+ \beta^2 \int_{t-t_0}^t \mu_s^{\otimes 2}\etc{P^{\otimes 2}_{t-s}f(\mom(s)) \un{\omega_1(s)=\omega_2(s)}}\, ds \notag\\
&-2 \beta^2 \int_{t-t_0}^t \mu_s^{\otimes 3}\etc{P^{\otimes 2}_{t-s}f(\mom(s)) (\un{\gamma(s)=\omega_1(s)} + \un{\gamma(s)=\omega_2(s)})}\, ds\notag \\
&+ 3 \beta^2 \int_{t-t_0}^t \mu_s^{\otimes 2}\etc{P^{\otimes 2}_{t-s}f(\mom(s))}\, I_s\, ds,\notag
\end{align}
where
$$ N_{t_0,t} =  \int_{t-t_0}^t \mu_s^{\otimes 2}\etc{P^{\otimes 2}_{t-s}f(\mom(s))(\beta \sum_i dB_{\omega_i(s)}(s) -2 \frac{dZ_s}{Z_s})}\,.
$$

The following inequalities are standard folklore,and are crucial in our proof: they will be used repeatedly hereafter and we provide a proof in the appendix.
\begin{equation}
  \label{ineq:r}
  0\le P^{\otimes 2}_t f(x_1,x_2) \le r(t) = P^{\otimes 2}_t f(x,x) \le 1
\end{equation}

In particular, we have 
\begin{align}
  \label{eqloc:un}
  I_t \ge  N_{t_0,t} & + \beta^2  \int_{t-t_0}^t r(t-s) I_s \, ds\\
& -4\beta^2 
 \int_{t-t_0}^t \mu_s^{\otimes 3}(P^{\otimes 2}_{t-s}f(\mom(s)) \un{\gamma(s)=\omega_1(s)})\, ds.\notag
\end{align}

Indeed, the second and fifth terms of~\eqref{eq:14} are non negative, in the second term we have

\begin{align*}
  P^{\otimes 2}_{t-s}f(\mom(s)) \un{\omega_1(s)=\omega_2(s)} &= P^{\otimes 2}_{t-s}f(\omega_1(s),\omega_1(s)) \un{\omega_1(s)=\omega_2(s)}\\
& = r(t-s) \un{\omega_1(s)=\omega_2(s)}\,,
\end{align*}

and finally, the fourth term can be written, thanks to symmetry of $f$,

$$-4\beta^2 
 \int_{t-t_0}^t \mu_s^{\otimes 3}(P^{\otimes 2}_{t-s}f(\mom(s)) \un{\gamma(s)=\omega_1(s)})\, ds\,.$$

{\bf Claim 1} : \begin{equation*}\mu_s^{\otimes 3}(P^{\otimes 2}_{t-s}f(\mom(s)) \un{\gamma(s)=\omega_1(s)}) \le I_s \inf(\sqrt{I_s r(t-s)}, r(t-s))\,.
\end{equation*}

\medskip
Indeed with $U_s(x)=\mu_s(\omega(s)=x)$ we have 
\begin{align*}
\mu_s^{\otimes 3}\etc{P^{\otimes 2}_{t-s}f(\mom(s)) \un{\gamma(s)=\omega_1(s)}} &=
\sum_x \mu_s^{\otimes 3}\etc{P^{\otimes 2}_{t-s}f(x,\omega_2(s))\, \un{\gamma(s)=\omega_1(s)=x}} \\
&= \sum_x U_s(x)^2 \mu_s(P^{\otimes 2}_{t-s}f(x,\omega(s)))
\end{align*}
and
\begin{equation*}
\mu_s(P^{\otimes 2}_{t-s}f(x,\omega(s))) = \sum_y U_s(y) P^{\otimes 2}_{t-s}f(x,y) \le r(t-s) \sum_y U_s(y) = r(t-s)\,.
\end{equation*}
We also have,  by Cauchy-Schwarz,
\begin{align*}
  \mu_s(P^{\otimes 2}_{t-s}f(x,\omega(s)))&\le \etp{\sum_y U_s(y)^2 \sum_y  (P^{\otimes 2}_{t-s}f(x,y))^2}^\undemi \\
&= \sqrt{I_s r(2(t-s))} \le \sqrt{I_s r(t-s)}\,,
\end{align*}
since if $\tilde{\omega}(t)= \omega_1(t)-\omega_2(t)$ we have, thanks to Markov property and symmetry, 
\newcommand{\tom}{\tilde{\omega}}

\begin{align*}
 r(2t) &= \prob{\tom(2t)=0} = \sum_y \PP_0(\tom(t)=y) \PP_y(\tom(t)=0) = \sum_y \PP_0(\tom(t)=y)^2 \\
&= \sum_y P^{\otimes 2}_tf(0,y)^2 = \sum_y P^{\otimes 2}_tf(x,y)^2
\end{align*}

{\bf Claim 2} :
\begin{align} \label{ito:3et}
4 \beta^2 R(t_0) \int_0^T J_s\, ds + \int_{t_0}^T I_s\,ds \ge& \int_{t_0}^T N_{t_0,t}\, dt \\
&+ 
\beta^2(1-4\sqrt{\epsilon_0})R(t_0) \int_{t_0}^{T- t_0} I_s\, ds\,.\notag
\end{align}

Observe that when $I_s \le c_0$ and $t-t_0 \le s \le t$, we have $I_s \le \epsilon_0 r(t-s)$, therefore, from Claim 1 we deduce that,
\begin{align*}
\int_{t-t_0}^t \mu_s^{\otimes 3}(P^{\otimes 2}_{t-s}f(\mom(s)) \un{\gamma(s)=\omega_1(s)})\, ds &\le \int_{t-t_0}^t {I_s \sqrt{I_s r(t-s)} \un{I_s \le c_0}}\, ds \\ &+ 
\int_{t-t_0}^t r(t-s) {I_s \un{I_s > c_0}}\, ds\,\\
&\le \sqrt{\epsilon_0} \int_{t-t_0}^t r(t-s)I_s\, ds\\& +  \int_{t-t_0}^t r(t-s)J_s\, ds.
\end{align*}
 Plugging this inequality into~\eqref{eqloc:un} yields

\begin{equation*}
  \label{eqloc:2}
I_t \ge N_{t_0,t} + \beta^2(1-4\sqrt{\epsilon_0})  \int_{t-t_0}^t r(t-s)I_s\, ds - 4 \beta^2  \int_{t-t_0}^t r(t-s)J_s\, ds\,.
\end{equation*}

Given $T\ge t_0$, we are going to integrate this inequality between $t_0$ and $T$. On the one hand,
\begin{align*}
\int_{t_0}^T dt \int_{t-t_0}^t r(t-s)J_s\, ds &=\int \int\un{0\le u\le t_0,t_0-u \le s \le T-u} J_s r(u)\,ds du \\
&\le R(t_0) \int_0^T J_s\, ds \,.
\end{align*}
On the other hand,
$$  \int_{t_0}^T dt \int_{t-t_0}^t r(t-s)I_s\, ds \ge \int_{t_0}^{T-t_0} I_s\, ds \int_0^{t_0} r(u)\, du = R(t_0) \int_{t_0}^{T-t_0} I_s\, ds \,.$$
The claim follows immediately.

\medskip

{\bf Claim 3} : let $\Nrond_T = \int_{t_0}^T N_{t_0,t}\, dt$. Then as ${T\to +\infty} $
$$ \frac{\Nrond_T}{\int_0^T I_s\, ds} \to 0 \quad \text{in probability}.$$

Let us defer the proof of this claim. Since $0\le I_s \le 1$ and $\intof I_s \, ds =+\infty$, we have, 
$$ \lim_{T\to +\infty} \frac{\int_{t_0}^T I_s\, ds }{\int_0^T I_s\, ds} = 
\lim_{T\to +\infty}  \frac{\int_{t_0}^{T-t_0} I_s\, ds}{\int_0^T I_s\, ds}=1\quad a.s.$$

Let $c_1 = \frac{\beta^2 (1 - 4 \sqrt{\epsilon_0})R(t_0) -1}{4 \beta^2 R(t_0)}$. If we divide \eqref{ito:3et} by $\phi_T=\int_0^T I_s\, ds$ and take $\limsup$ as $T\to +\infty$, we obtain that \as
\begin{align*}
 \limsup_{T\to \infty} \unsur{\phi_T} \int_0^T J_s\, ds -c_1& \ge
\limsup_{T\to \infty}  \frac{\Nrond_T}{ 4 \beta^2 R(t_0)\phi_T} \\
& \ge \limsup_{T\to +\infty} - \frac{\valabs{\Nrond_T}}{4 \beta^2 R(t_0)\phi_T} \\
&= - \liminf_{T\to +\infty}  \frac{\valabs{\Nrond_T}}{4 \beta^2 R(t_0)\phi_T}\\
&=0\,.
\end{align*}

This yields

$$  \limsup_{T\to \infty}\frac{\int_0^T J_s ds}{\int_0^T I_s\, ds} \ge c_1 \quad a.s. $$

\bigskip

{\bf Proof of Claim 3.}

 By Fubini's theorem,
\begin{eqnarray*}
{\Nrond}_T &=& \int_{t_0}^T dt \int_{t-t_0}^t \mu_s^{\otimes 2}\left[ P^{\otimes 2}_{t-s} f(\omega_1(s), \omega_2(s))
\big( \sum_i \beta d B_{\omega_i(s)} (s) - 2 \frac{ d Z_s}{ Z_s}\big) \right] \\
   &=& \int_0^T \mu_s^{\otimes 2}\left[ G(s,\omega_1(s), \omega_2(s))
\big( \sum_i \beta d B_{\omega_i(s)} (s) - 2 \frac{ d Z_s}{ Z_s}\big) \right],
\end{eqnarray*}

\noindent with $$ 0 \le G(s, x_1, x_2) :=  \int_{(t_0- s)^+}^{(T-s)^+ \wedge t_0} P^{{\otimes 2}}_{t-s}f(x_1, x_2)\, dt
\le t_0, \quad \forall x_1, x_2 \in {\mathbb Z}^d.$$

Let us view   $\Nrond_T=X_T$ as the value at time $T$ of the continuous
martingale
 
$$X_t=\int_0^t  \mu_s^{\otimes 2}\left[ G(s,\omega_1(s), \omega_2(s))
\big( \sum_{i=1}^2 \beta d B_{\omega_i(s)} (s) - 2 \frac{d Z_s}{Z_s}\big) \right]\,.$$

 We can compute its quadratic variation :

$$ \crochet{X,X}_T  \le 4 \beta^2 \int_0^T  \mu_s^{\otimes 4} \left [    G(s,
\omega_1(s) , \omega_2(s))  G(s, \omega_3(s) , \omega_4(s)) \big(   1_{(\omega_1(s)=\omega_3(s))} +
I_s \big) \right] d s ,$$

\noindent which satisfies  
\begin{equation}
  \label{eq:13}
\crochet{X,X}_T  \le 8 \beta^2 t_0^2 \int_0^T  I_s ds .
\end{equation}

 Let $\epsilon>0$, we shall prove that  \begin{equation}
 \lim_{T \to \infty}  {\bf P}\Big( \Nrond_T > \epsilon \, \int_0^T I_s d s
 \Big) =0. \label{cl3}
 \end{equation}

 To this end, define $\delta = \epsilon/(8\beta^2 t_0)$. 
  We have   $$ \gesp{e^{\delta \Nrond_T -\frac{\delta^2}{ 2}\crochet{X,X}_T}}= {\bf E}\etc{ e^{
 \delta X_T - \frac{\delta^2}{ 2} \langle X,
 X\rangle_T}}=1.$$

(since $\crochet{X,X}_T$ is bounded, Novikov's criterion for the exponential martingale is obviously satisfied). It follows that 

\begin{align*} 1 &\ge {\bf E} \Big(
1_{(\Nrond_T >\epsilon \, \int_0^T I_s d s)} e^{ \delta \Nrond_T -
\frac{\delta^2}{ 2} \langle X,
 X\rangle_T}\Big) \\
    &\ge {\bf E} \Big(
1_{(\Nrond_T >\epsilon \, \int_0^T I_s d s)} e^{ ( \delta \epsilon
  - \frac{\delta^2}{ 2} 8\beta^2 t_0)  \int_0^T I_s d s}\Big) \\
    &= {\bf E} \Big(
1_{(\Nrond_T >\epsilon \, \int_0^T I_s d s)} e^{   4 \beta^2 t_0
\delta^2 \int_0^T I_s d s}\Big) &\text{by \eqref{eq:13}}\\
    &\ge  e^{ 4 \beta^2 t_0 \delta^2  K} \, {\bf P}\Big( \Nrond_T >\epsilon \, \int_0^T I_s ds, \, \int_0^T I_s d     s \ge K\Big),
\end{align*}

for any constant $K>0$.  Consequently, we have
$${\bf P}\Big( \Nrond_T >\epsilon \, \int_0^T I_s d    s \Big) \le
 {\bf P}\Big(  \int_0^T I_s d    s < K  \Big) + e^{ -4 \beta^2 t_0 \delta^2
 K}.$$

Since $\int_0^T I_s d    s \to \infty$ almost surely, we get
$$ \limsup_{T\to\infty} {\bf P}\Big( \Nrond_T >\epsilon \, \int_0^T I_s d    s
\Big)\le e^{ -4 \beta^2 t_0 \delta^2
 K},$$

for any constant $K>0$.  Then    by letting $K\to\infty$ we get
(\ref{cl3}). Considering the martingale $-X$, we prove in the same way that
\begin{equation}
 \lim_{T \to \infty}  {\bf P}\Big(- \Nrond_T > \epsilon \, \int_0^T I_s d s
 \Big) =0. 
 \end{equation}

 and this complete the proof of {\bf Claim 3}.

\section*{Appendix}

We provide a proof of~\eqref{ineq:r}. Recall that $f(x,y)=\un{x=y}$. We let $p_t(x)=\prob{\omega(t)=x}$ be the distribution of simple random walk at time $t$. Then, by translation invariance:
\begin{align*}
   P^{\otimes 2}_t f(x_1,x_2)&= \PP^{\otimes 2}_{x_1,x_2}\etp{\omega_1(t)
=\omega_2(t)} \\
&= \PP^{\otimes 2}\etp{x_1+\omega_1(t)=x_2+\omega_2(t)}\\
&=\sum_z \prob{x_1 + \omega_1(t)=z} \prob{x_2 +\omega_2(t)=z} &\text{(by independence)}\\
&= \sum_z p_t(z-x_1)p_t(z-x_2) \\
&\le \etp{\sum_z p_t(z-x_1)^2}^\undemi \etp{\sum_z p_t(z-x_2)^2}^\undemi &\text{(by Cauchy-Schwarz)}\\
&= \sum_z p_t(z)^2 = r(t)\,.
\end{align*}

\section*{Acknowledgements} We want to thank both  referees for a careful reading of the first version of this paper and for suggestions for improvements for the exposition. We are also grateful to A. Camanes for his comments on a very early draft of this paper.

\bibliographystyle{siam}
\bibliography{new,poly}

\end{document}